\documentclass[10pt,reqno]{amsart}
\usepackage{latexsym}
\usepackage{graphicx}
\usepackage{fancyhdr,amssymb}
\usepackage{amsmath}
\usepackage{tikz}
\usepackage{pgf}
\usetikzlibrary{automata, positioning}

\usepackage{hyperref} 

\newcommand{\N}{{\mathbb N}}

\theoremstyle{plain}
\numberwithin{equation}{section}
\newtheorem{thm}{Theorem}[section]
\newtheorem{theorem}[thm]{Theorem}
\newtheorem{lemma}[thm]{Lemma}

\theoremstyle{definition}

\newtheorem{conjecture}[thm]{Conjecture}  % additional proclamation for this paper

\begin{document}
\fancyhead{}
\renewcommand{\headrulewidth}{0pt}
\fancyfoot{}
\fancyfoot[LE,RO]{\medskip \thepage}
\fancyfoot[LO]{\medskip MISSOURI J.~OF MATH.~SCI., FALL 2022}
\fancyfoot[RE]{\medskip MISSOURI J.~OF MATH.~SCI., VOL.~34, NO.~2}

\setcounter{page}{1}

\title{An Abelian Loop for Non-Composite Numbers}
\author{Raghavendra N Bhat}
\address{\newline Raghavendra N Bhat\newline 
University of Illinois, Urbana Champaign \newline
Department of Mathematics\newline
1409 West Green Street\newline
Urbana, IL 61801
}
\email{rnbhat2@illinois.edu}

\begin{abstract}
We define an abelian loop on a set \(S\) consisting of 1 and all odd prime numbers with an operation $\bullet$, where for $a,b$ $\in$ $S$, $a$  $ \bullet$ $b$ is the smallest element of \(S\)
strictly larger than $|a-b|$. We use theorems and conjectures from number theory to prove properties of the loop and state analogous conjectures about the loop.
\end{abstract}

\maketitle

%-------------------------------------------------------

\section{Introduction}

%-------------------------------------------------------
Let \(S\) consist of 1 and all odd prime numbers, with operation $\bullet$, where for $a,b$ $\in$ $S$, $a$  $ \bullet$ $b$ gives the smallest element of \(S\) strictly larger than $|a-b|$.
This set is an abelian loop, because it has the following properties:
\begin{enumerate}
\item $S$ is closed under $\bullet$, since the operation always generates non-composite numbers in $S$.
\item $\bullet$ is commutative, owing to the absolute value in the operation. Thus, $a$ $\bullet$ $b = b$ $\bullet$ $a$
\item 1 is the identity element of $S$ since for all $a$ $\in$ $S$, $a$ $\bullet$ 1 = $a$.
\item $S$ is also closed under inversion because for all $a$ $\in$ $S, a$ $\bullet$ $a = 1$. Hence, the inverse of an element in $S$ is the element itself.
\end{enumerate}
Thus, we have an abelian loop consisting only of 1 and all odd primes. $S$ is not a group because there is no associativity. For simplicity, we will define a function $N(x)$ to be the smallest element of $S$ that is greater than $|x|$.

%-------------------------------------------------------

\section{Some Results}

%-------------------------------------------------------

We now apply some well known results in number theory to $S$ to prove some theorems.

\begin{theorem}
For all $s \in S,$ there exists $a \in S$ such that $a \bullet s = a.$
\end{theorem}

\ \\
\noindent\rule{0.84in}{0.4pt} \par
\medskip
\indent\indent {\fontsize{8pt}{9pt} \selectfont DOI: 10.35834/YYYY/VVNNPPP \par}
\indent\indent {\fontsize{8pt}{9pt} \selectfont MSC2020: 11P32, 20N05 \par}
\indent\indent {\fontsize{8pt}{9pt} \selectfont Key words and phrases: loops, differences of primes. \par}

\thispagestyle{fancy}

\vfil\eject
\fancyhead{}
\fancyhead[CO]{\hfill AN ABELIAN LOOP FOR NON-COMPOSITE NUMBERS}
\fancyhead[CE]{R.~N.~BHAT  \hfill}
\renewcommand{\headrulewidth}{0pt}

\begin{proof}
It is well known that there are arbitrarily long strings of composite numbers. Thus, there is a prime $a$ that is preceded by at least $2s$ composites. Then $a$ $\bullet$ $s$ = $N(a - s)$ = $a$,
because $a - s$ is in the string of $2s$ composites.
\end{proof}

\begin{theorem} \label{T:chain}
For every $n$ $\in$ $\N$ there exist $a_{1}, a_{2}, a_{3}, \dots, a_{n}$ $\in S$, all distinct, such that $a_{1} \bullet a_{2} = a_{2} \bullet a_{3} = a_{3} \bullet a_{4}= \dots = a_{n-1} \bullet a_{n}$.
\end{theorem}
\begin{proof}
The theorem is a direct derivation from Green-Tao's theorem \cite{Green} of 2004. Our aim is to create an arbitrarily long chain of non-composites such that each adjacent pair has the same value when the $\bullet$ operation is performed.

Green and Tao proved that for any natural number $n$, there exists an arithmetic progression of $n$ primes. Since all primes except 2 are present in our group $S$, we can make use of the theorem. Let $a_{1}, a_{2}, a_{3}, \dots, a_{n}$ be an arithmetic progression of $n$ primes. Let their common difference be $k$. Let $l$ be the smallest non-composite number larger than $k$. Thus, for every pair $(a_{j}, a_{j+1})$ from our arithmetic progression, we have $a_{j} \bullet a_{j+1} = l.$ The proof is complete.
\end{proof}

It is important to note that the number 2 is not part of any arithmetic progression of primes for $n \geq 3$ owing to parity issues.

Here is an example of Theorem \ref{T:chain}. For $n=4,$ we have the arithmetic progression: 41, 47, 53, 59. 
Then 41 $\bullet$ 47 = 47 $\bullet$ 53 = 53 $\bullet$ 59 = 7.

\begin{theorem} \label{T:notriangle}
It is not possible to have $a, b, c$ $\in$ $S$ such that $a$ $\bullet$ $b = b$ $\bullet$ c = a $\bullet$ $c$, where $a, b, c$ are different from one another.
\end{theorem}

\begin{proof}
Assume, for the sake of contradiction, that $a \bullet b = b \bullet c = a \bullet c.$ Without loss of generality, let $a<b<c$. Thus, the smallest prime strictly greater than $|a-b|$, $|b-c|$ and $|a-c|$ should be the same number. Assume this to be some prime $p$. Let $|a-b| = x$, $|b-c| = y$ and $|a-c| = x+y$. Thus, we have the nearest prime to $x, y, $and $x+y$ equal $p$.

This implies that there is no prime between $x$ and $x+y$, and between $y$ and $x+y$.
But $x+y$ is $\geq$ 2(min$\{x,y\})$. This contradicts Bertrand's Postulate \cite{Dickson} which states that for every number $n$ $\geq$ 2, there exists a prime between $n$ and $2n$. We arrive at a contradiction.
Thus, it is not possible to have $a \bullet b = b \bullet c = a \bullet c$ for while $a, b, c$ are different from one another.
\end{proof}

\begin{lemma} \label{L:adjacent}
If $t$ is a positive even integer, then $N(t)$ and $N(t + 2)$ are either equal or are consecutive in $S$.
\end{lemma}

\begin{proof}Let $s_1$ = $N(t)$ and $s_2$ be the next element in $S$ after $s_1$. Since all the elements of $S$ are odd, $s_2 \geq s_1 + 2 = N(t) + 2 \geq t + 2$. Thus, $s_1 = N(t) \leq N(t+2) \leq s_2$. Since $s_1$ and $s_2$ are consecutive, $N(t+2)$ is one of them.
\end{proof}

\begin{theorem} \label{T:twinprimes}
If $a$ and $b$ are twin primes in $S$, then  $\forall x$ $\in$ $S$ such that $x < a$ and $x < b$, either
\begin{enumerate}
\item
$a \bullet x = b \bullet x$, or
\item
$a \bullet x$ and $b \bullet x$ are adjacent in $S$.
\end{enumerate}
\end{theorem}

\begin{proof}
Let $a$ and $b$ be twin primes, with $a < b$, without loss of generality. Let $x$ be an element of $S$ smaller than both. If we let $t=a-x$ and $t+2=b-x$, by our Lemma \ref{L:adjacent}, we have $N(t)$ and $N(t+2)$ to be equal or consecutive in $S$. The proof is complete.
\end{proof}

%-------------------------------------------------------

\section{Conjectures}

%-------------------------------------------------------

We now take a present a few conjectures in our loop $S$.

\begin{conjecture} \label{C:likeTwinPrime}
$a \bullet b = x$ has infinitely many solutions for all $x$ in $S$.
\end{conjecture}

This is a straight forward analogy of the Polignac conjecture \cite{Polignac} which claims infinite prime pairs for all even gaps. For any $x$ in $S$, pick even number $k$, such that $N(k) =x$. This establishes the equivalence of Conjecture \ref{C:likeTwinPrime} to the Polignac conjecture.

We can restate the twin prime conjecture as the existence of infinitely many solutions to $a \bullet b = 3$. The equivalence can be proved as follows.

Suppose $p_1 < p_2$ are odd primes. If
they are twin primes, then then $N(p_2 - p_1) = N(2) = 3$. If $N(p_2 - p_1) = 3$, then $p_2 - p_1 < N(p_2 - p_1) = 3$, thus $p_2 - p_1 = 2$.

\begin{conjecture} \label{C:quartet}
The quartet $a, b, c, d$ (distinct numbers) such that $a \bullet b = b \bullet c = c \bullet d = d \bullet a$ has infinitely many solutions.
\end{conjecture}

Here are a few examples:
\begin{enumerate}
\item
101, 7, 103, 11. We have 101 $\bullet$ 7 = 7 $\bullet$ 103 = 103 $\bullet$ 11 = 11 $\bullet$ 101 = 97.
\item
137, 23, 139, 19. We have 137 $\bullet$ 23 = 23 $\bullet$ 139 = 139 $\bullet$ 19 = 19 $\bullet$ 137 = 127.
\end{enumerate}
The intuition for the conjecture comes from both the twin prime conjecture and the Polignac conjecture for $n=4$. Primes following the latter are sometimes referred to as cousin primes. An example would be the pair $(19,23)$.

Assume $a$ and $c$ to be twin primes. From Theorem \ref{T:twinprimes}, we know that it is possible to have $b$ $\in S$ such that $a \bullet b$ equals $b \bullet c$. The conjecture hopes for the existence of $d$ such that $a \bullet d$ equals $d \bullet c$. In addition, we hope that $b$ and $d$ are not too far from each other, hence allowing the possibility of $a \bullet b$ equalling $a \bullet d$. Letting $b$ and $d$ be twins or cousins increases the chance of finding examples for the conjecture. Owing to little knowledge about the infinite existence of twins or cousin primes, we will explore these ideas visually in the next section.

%-------------------------------------------------------

\section{Geometrical Analysis}

%-------------------------------------------------------

To elaborate Conjecture \ref{C:quartet}, we now use shapes to introduce a new way of visualizing the operation $\bullet$ in the loop S. The set up is as follows:
\begin{enumerate}
\item Each point on the plane defines an element of $S$.
\item A line connecting two points on the plane will have a length (distance) that is the value of the $\bullet$ operation between the two points.
\item We will draw shapes such that each vertex of the shape is labelled with an element of $S$ and each side has length $A \bullet B$ where $A$ and $B$ are the labels of its endpoints.
\end{enumerate}
We now address the following questions : Which $n$-sided polygons can be created such that all their sides have the same value? As proved in Theorem \ref{T:notriangle}, an equilateral triangle is not possible.

We will use the example from Conjecture \ref{C:likeTwinPrime} to demonstrate a rhombus.

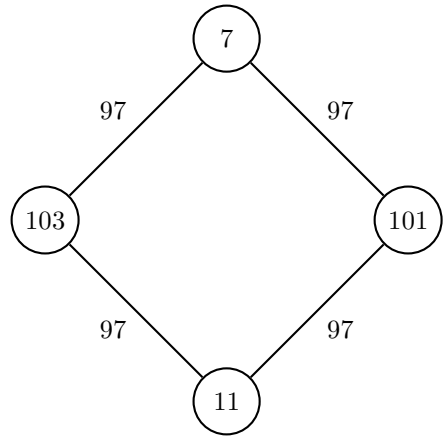
\begin{center}
\begin{tikzpicture}[auto,node distance=2.5cm,
   thick,main node/.style={draw,font=\bfseries}]
    \node[state] (7) {$7$};
    \node[state] (101) [below right=of 7]  {$101$};
    \node[state] (103) [below left =of 7]  {$103$};
    \node[state] (11) [below right =of 103]  {$11$};
    \path[-]
    (7) edge  node {97} (101)
    (103) edge   node {97} (7)
    (11) edge  node {97} (103)
    (101) edge  node {97} (11);
\end{tikzpicture}
\end{center}
From Theorem \ref{T:notriangle}, we also infer that a rhombus with a diagonals having the same length as it sides is not possible as it creates two equilateral triangles.
The Green-Tao analogy used in the example to Theorem \ref{T:chain} can be represented as follows:

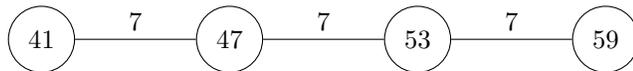
\begin{center}
\begin{tikzpicture}[shorten >=0.005pt,node distance=2.5cm,on grid,auto]
    \node[state] (41) {$41$};
    \node[state] (47) [right=of 7]  {$47$};
    \node[state] (53) [right=of 47]  {$53$};
    \node[state] (59) [right=of 53]  {$59$};
    \path[-]
    (41) edge  node {7} (47)
    (47) edge  node {7} (53)
    (53) edge  node {7} (59);
    
\end{tikzpicture}
\end{center}
It is important to note here that the primes need not necessarily be in an arithmetic progression. All we require is the nearest prime to the difference of two adjacent vertices to be the same.

\medskip
\noindent
{\textbf{Acknowledgements.}}
I would like to acknowledge the professors, at the University of Illinois for their valuable advice and support.

\end{document}